\documentclass[journal]{IEEEtran}

\usepackage{amsmath,}
\usepackage{array}
\usepackage{amssymb}
\usepackage{amsfonts}
\usepackage[]{algorithm2e}
\usepackage[font=small,labelfont=bf]{caption}
\usepackage{multirow}
\usepackage{cite}
\usepackage{textcomp}
\usepackage{float}
\usepackage{epsfig}
\usepackage{epstopdf}
\usepackage{color}
\usepackage{amsthm}
\usepackage{booktabs}
\usepackage{xspace}
\usepackage{xspace,amssymb,epsfig,syntonly}
\usepackage{epsfig,amsmath,color}
\usepackage{xspace,syntonly}%,empheq}
\usepackage{cite}
\usepackage{graphicx}
\usepackage{soul}
\usepackage{hyperref}

\newtheorem{theorem}{Theorem}

\newtheorem{remark}[theorem]{Remark}
\newtheorem{claim}[theorem]{Claim}

\DeclareMathOperator{\real}{Re}
\DeclareMathOperator{\imag}{Im}

\title{Beyond Relaxation and Newton-Raphson: Solving AC OPF for Multi-phase Systems with  Renewables}
\author{Ahmed S. Zamzam%, {\it Student Member, IEEE}
	, Nicholas D. Sidiropoulos%, {\it Fellow, IEEE}
	, and Emiliano Dall'Anese%, {\it Member, IEEE}
\thanks{A. S. Zamzam and N. D. Sidiropoulos are with the Department of Electrical and Computer Engineering, Univ. of Minnesota, Minneapolis, MN 55455, USA. Emails: {\small \tt\{zamza002,nikos\}@umn.edu}}%
	\thanks{E. Dall'Anese is with the National Renewable Energy Laboratory, Golden, CO 80401, USA. Email: {\tt\small emiliano.dallanese@nrel.gov}}%
	\thanks{The work of A.S. Zamzam and N.D. Sidiropoulos was partially supported by NSF under grants ECCS-1231504 and CIF-1525194. The work of E. Dall'Anese was supported in part by the Laboratory Directed Research and Development Program at the National Renewable Energy Laboratory.}
}
\begin{document}
\maketitle
\begin{abstract}

This paper focuses on the AC \textit{Optimal Power Flow} (OPF) problem for  multi-phase systems. Particular emphasis is given to systems with high integration of renewables, where adjustments of the real and reactive output powers from renewable sources of energy are necessary in order to enforce voltage regulation.  The AC OPF problem is known to be nonconvex (and, in fact, NP-hard). Convex relaxation techniques have been recently explored to solve the OPF task with reduced computational burden; however, sufficient conditions for tightness of these relaxations are only available for restricted classes of system topologies and problem setups. Identifying {\em feasible} power-flow solutions remains hard in more general problem formulations, especially in unbalanced  multi-phase systems with renewables. To identify feasible and optimal AC OPF solutions in challenging scenarios where existing methods may fail, this paper leverages the Feasible Point Pursuit - Successive Convex Approximation algorithm -- a powerful approach for general nonconvex quadratically constrained quadratic programs. The merits of the approach are illustrated using single- and multi-phase distribution networks with renewables, as well as several transmission systems.
\end{abstract}

\begin{keywords}
Optimal power flow, renewable sources of energy, convex relaxation, feasible point pursuit, successive convex approximation.
\end{keywords}

%\vfill

\section{Introduction}
The AC optimal power flow (OPF) problem is a predominant task in optimizing the performance of power grids. The OPF problem aims at minimizing an appropriate operational cost while respecting the network's physical and engineering constraints. However, due to the quadratic nature of the power flow equations, the OPF problem is known to be nonconvex and NP-hard in general~\cite{Lehmann2016, Lavaei-2012}. Existing approaches to the OPF problem range from classical Newton-Raphson, to Lagrangian relaxation, genetic algorithms and interior point methods. Unfortunately, these methods do not  provide optimality or feasibility guarantees except in certain cases, and are quite sensitive to the initial guess. {While the Newton-Raphson method has been traditionally employed to solve AC OPF problems for transmission systems, its convergence is challenged when it is applied to multi-phase distribution networks; this is primarily due to the high resistance-to-reactance ratio of distribution lines, which can cause the Jacobian matrix to be ill-conditioned. }

Many recent research efforts have been trying to approach the solution of OPF problem using relaxation techniques~\cite{Bai-2008, Lavaei-2012, Low-2014-I, Jabr-2006, Farivar-2013, Coffrin2015qc, Zhang2013relaxed, Gan2014convex, Bienstock2015lp, Zhang15}. Among those, semidefinite relaxation (SDR) was shown to be able to find the global optimal solution of the problem in many cases. SDR relies on matrix-lifting and rank relaxation to convexify the feasible set of the OPF problem~\cite{Bai-2008, Lavaei-2012, Louca-2013, DallAnese13, Zhang15, Madani15}; the resulting relaxed problem can be solved in polynomial time. OPF-optimality of the SDR solution can always be tested {\em a posteriori} by checking the rank of the SDR solution matrix; but it is very useful to know {\em a priori} in which cases SDR will yield an optimal solution for the original nonconvex OPF problem. These are the cases when SDR yields a solution that is rank-one, or can be easily transformed to rank-one. In those cases, SDR is not a relaxation after all; we say that SDR is {\em tight}. Tightness of SDR relaxation was proved for a number of network setups under restrictive conditions. In~\cite{Lavaei-2012}, tightness of SDR was shown for a network comprising only resistive loads, provided load over-satisfaction is allowed and the dual variables are all positive. Assuming availability of sufficient phase shifters, it was proven that SDR is exact if load over-satisfaction is allowed~\cite{Sojoudi-2012}. For tree/radial networks, under operational constraints on voltage magnitudes, line losses, and line flows, the SDP relaxation was shown to be tight if there are no lower limits on the power generation~\cite{Zhang15}. This result was extended in~\cite{Lavaei-2014} for radial networks with lower limits only on the active power, under reasonable conditions. The inexactness of SDR for a general network was demonstrated in~\cite{Lesieutre-2011} using a simple $3$-bus network. Unfortunately, in cases where SDR is not tight, it is not easy to recover a physically meaningful solution from the solution matrix; only a lower bound on the optimal cost is provided. An approach to recover an OPF solution from the SDR solution was pursued  in~\cite{Louca-2013}, but still there is no guaranty of recovering a physically meaningful OPF solution.

Another relaxation technique was proposed in~\cite{Jabr-2006} for radial networks. The method eliminates the voltage angles by defining new variables representing the real and imaginary parts of the second order voltage moments, and then expresses the power flow equation in terms of the new variables. \cite{Jabr-2006} solves the problem using second order cone programming (SOCP). Tightness of this SOCP relaxation is an open issue. Along the same lines, a quadratic convex (QC) relaxation approach was proposed in~\cite{Coffrin2015qc} where the network constraints are replaced by convex surrogates. Although shown to be tighter than the SOCP relaxation, the QC relaxation also lacks proof of tightness, and can return solutions that are infeasible for the original OPF problem.

As a generalization of SDR, moment-based relaxation has been proposed in~\cite{Molzahn-2014I} using the Lassarre Hierarchy. Empirically, the method has been demonstrated to be able to find the OPF solution in cases where SDR fails. The moment-based relaxation introduces higher order voltage moments as new variables and defines the power flow quantities in terms of these moments. On the down side, considering higher order moments requires solving very large SDP instances which may not be computationally tractable. Aiming to alleviate the computational burden of moment-based relaxation,~\cite{Molzahn-2015I} exploited the structure of the OPF problem to develop a more tractable solution for low order moments. However, due to the NP-hardness of the problem, the moment order required to approach the optimal solution may be very large. Building on the same tool, a Laplacian-based approach has been proposed in~\cite{Molzahn-2016I}, where an upper bound on the cost function is assumed, and the cost is replaced by a function that penalizes constraint violations.

It is also worth emphasizing that the AC OPF task is becoming increasingly important for distribution systems with high integration of renewable energy resources (RESs), where adjustments of the real and reactive output-powers from renewable sources are necessary in order to enforce voltage regulation. Particularly relevant is the case of distribution feeders with high penetration of photovoltaic (PV) systems, where reverse power flows induced by PV-systems operating according to current practices may increase the likelihood of overvoltage conditions.  OPF formulations aim  at minimizing the cost of real power curtailment as well as the cost for reactive power support, while concurrently pursuing utility-oriented objectives and  ensuring voltage regulation~\cite{Dallanesse-2014, Su14,Bonfiglio14,Guggilam16,Wang-2016}. However, in this particular setting,  the overall cost function of the OPF task may not be strictly increasing in the power injections, which implies that relaxation methods such as SDR~\cite{DallAnese13,Gan14,Robbins16} are not guaranteed to be tight.

%\textcolor{blue}
{Overall}, the main contributions of this paper are as follows:
\begin{enumerate}
\item {\em Designing an efficient algorithm that can solve the OPF problem when the relaxation approaches fail to find an optimal solution.} Inspired by recent advances in solving nonconvex quadratically-constrained quadratic programs (QCQP), the OPF problem is formulated as a nonconvex QCQP. In~\cite{Mehanna-2014}, a \textit{Feasible Point Pursuit Successive Convex Approximation} (FPP-SCA) algorithm was proposed, and empirically shown to be very effective in solving nonconvex QCQP problems in cases where SDR fails. The FPP-SCA algorithm replaces the nonconvex constraints by inner convex surrogates around a specific point to construct a convex restriction of the original problem. Such restriction may lead to infeasibility, even if the original problem is feasible. The main idea behind FPP-SCA is to allow a controllable amount of constraint violations to enable the algorithm to make progress towards feasibility in its initial stages. Towards this end, a slack variable is added to ensure feasibility at each step, and the cost function is augmented with a term that penalizes the slack that reflects the constraints violations. The overall approach is neither restriction nor relaxation, but rather judicious {\em approximation} of the OPF problem in each iteration, the solution of which is subsequently used as the approximation point for the next iteration. Upon finding a feasible voltage profile, successive convex approximation of the feasible set is used to find a KKT point of the OPF problem.
\item { {\em Identifying OPF solutions when minimizing the cost of active power curtailment and reactive power support.} The modified problem is solved to obtain an optimal voltage profile that conforms to the power system's operational and economic constraints. SDR is very sensitive to the choice of the cost function (especially when the cost function is non-increasing with the power flows in the network). On the other hand, the proposed algorithm is shown to be an effective approach for solving the modified OPF problem for single-phase and multi-phase system models.}
\item {\emph{Performance comparison}. The performance of the FPP-SCA is benchmarked against existing convex relaxation approaches and the IPOPT solver. Results demonstrate that the FPP-SCA algorithm is able to find solutions that are optimal or near optimal, even in cases where convex relaxation approaches and the IPOPT solver fail. }
{\item {\em Identifying problematic constraints when the AC OPF problem is infeasible.} The FPP algorithm has the ability to identify problematic constraints in cases where the AC OPF problem is infeasible. This is a distinct feature of the proposed method that off-the-shelf solvers such as IPOPT do not offer. }
\end{enumerate}

%A conference version of this paper has been submitted for publication to IEEE CDC 2016 \cite{Zamzam-2016}. The conference version only dealt with single-phase transmission systems without considering renewables. Only the cost of power generation was considered. Relative to \cite{Zamzam-2016}, this paper includes modeling and optimization for multi-phase systems. Incorporating renewable energy resources in the power network is also included. This work also accounts for the cost associated with active power curtailment and reactive power injections.

\noindent {\bf Organization:} The rest of the paper is structured as follows. The AC OPF problem formulation is introduced in Section II for multi-phase network with renewables, where the problem is reformulated as nonconvex QCQP. Section III describes the application of the FPP-SCA algorithm to the AC OPF problem in two phases.
%{followed by a modified version of the algorithm to identify problematic constraints}.
Test cases using different three-phase and single-phase networks are used to show the efficacy of the proposed algorithm in Section IV. Conclusions are drawn in Section V. %\reminder{From Nikos: add some more detail to the organization?}

\noindent {\bf Notation:} matrices (vectors) are denoted by boldface capital (small) letters; $(\cdot)^T$, $ (\cdot)^\ast $ and $(\cdot)^H$ stand for transpose, conjugate and complex-conjugate transpose, respectively; and $|(\cdot)|$ denotes the magnitude of a number or the cardinality of a set.

\section{Problem Formulation}

Consider a multi-phase network comprising $ N + 1 $ buses. The system is modeled by a graph $ \mathcal{G}:=(\mathcal{N}, \mathcal{L}) $, where $ \mathcal{N}:= \{0, 1, 2, \cdots, N\} $ is the set of multi-phase buses (nodes) and $ \mathcal{L} \subseteq \mathcal{N} \times \mathcal{N}$ represents the set of lines. Let bus 0 be the reference bus, whose voltages are taken as a reference for the phasorial representation. The set of phases of node $ k $ and phases of line $ (l, m) $ are denoted by $ \boldsymbol{\varphi}_k $ and $ \boldsymbol{\varphi}_{lm} $, respectively. Let $ v_{k,\phi} \in \mathbb{C}$ and $ i_{k,\phi} \in \mathbb{C} $ denote the phasor for the line-to-ground voltage and the current at node $ k $ for phase $ \phi $, and define $ \mathbf{v}_k :=[v_{k,\phi}]_{\phi \in \boldsymbol{\varphi}_k} $ and $ \mathbf{i}_k :=[i_{k,\phi}]_{\phi \in \boldsymbol{\varphi}_k} $. For notational simplicity, the paper hereafter focuses on three-phase systems; however, the proposed framework is applicable to systems featuring a variety of three-, two-, and single-phase nodes and branches.

Conventional fossil-fuel generators are assumed to be located at nodes $ \mathcal{G} \subseteq \mathcal{N} $, with $ P_{k,\phi}^{(G)} $, $ Q_{k,\phi}^{(G)} $ denoting the active and reactive power generated at phase $ \phi $ of bus $ k \in \mathcal{G}$. The load connected to phase $ \phi $ at bus $ k $ is denoted by $ P_{k,\phi}^{(L)} + j Q_{k,\phi}^{(L)} \in \mathbb{C} $. In addition, the apparent power transferred from bus $ l \in \mathcal{N} $ to the rest of the network through line $ (l,m) \in \mathcal{L}$ for phase $ \phi $ is given by $ S_{lm,\phi} = P_{lm,\phi} + jQ_{lm,\phi}$.

Subset $ \mathcal{R} \subset \mathcal{N} $ collects nodes with installed renewable energy sources (RESs) such as PV systems. Given prevailing ambient conditions, let the available active power from the RES located at phase $ \phi $ of bus $ k \in \mathcal{R}$ be denoted by $ \overline{P}_{k,\phi}^{(R)} $. Also, let $ P_{k,\phi}^{(R)} $ and $ Q_{k,\phi}^{(R)} $ denote the injected active power and the injected/absorbed reactive power at bus $ k $ for phase $ \phi $. It is assumed that both active and reactive output-powers are controllable~\cite{Dallanesse-2014,Su14}. Accordingly, the allowed operating region of an RES can be described as follows:
\begin{equation}
\mathbf{\Psi}_{k,\phi} := \left \{P_{k,\phi}^{(R)}, Q_{k,\phi}^{(R)} : \begin{array}{c}
0\leq P_{k,\phi}^{(R)} \leq \overline{P}_{k,\phi}^{(R)}\\
(P_{k,\phi}^{(R)})^2 + (Q_{k,\phi}^{(R)})^2 \leq \overline{S}_{k,\phi}^2\\
|Q_{k,\phi}^{(R)}| \leq \tan (\overline{\theta}_{k,\phi}) P_{k,\phi}^{(R)}
\end{array}\right\}
\end{equation}
where $ \overline{S}_{k,\phi}^{(R)} $ represents the RES-inverter capacity, and $ \overline{\theta}_{k,\phi} $ capture minimum power factor requirements.

%Define $ \boldsymbol{\varphi}_k \equiv \boldsymbol{\varphi}_{lm} := \{a, b, c\}$ for all bus $ k \in \mathcal{N} $ and line $ (l, m) \in \mathcal{L}$, and hence $ \mathbf{v}_k := [v_{k,a}\quad v_{k,b}\quad v_{k,c}] $ and $ \mathbf{i}_k := [i_{k,a}\quad i_{k,b}\quad i_{k,c}] $.
Collect  voltages and currents in the  vectors $ \mathbf{v}:= [\mathbf{v}_0^T, \mathbf{v}_1^T, \cdots, \mathbf{v}_N^T]^T $ and $ \mathbf{i}:= [\mathbf{i}_0^T, \mathbf{i}_1^T, \cdots, \mathbf{i}_N^T]^T $ of length $ 3(N+1) $, respectively. Lines $ (l,m) \in \mathcal{L} $ are modeled as $ \pi $-equivalent circuit, where the phase impedance and shunt admittance matrices are denoted by $ \mathbf{Z}_{lm} \in \mathbb{C}^{|\boldsymbol{\varphi}_{lm}| \times |\boldsymbol{\varphi}_{lm}|} $ and $ \mathbf{\overline{Y}}_{lm} \in \mathbb{C}^{|\boldsymbol{\varphi}_{lm}| \times \boldsymbol{\varphi}_{lm}} $, respectively. Voltages and injected currents abide by Ohm's law and Kirchhoff's law, which lead to the compact relationship $ \mathbf{i} = \mathbf{Y} \mathbf{v}$. The network admittance matrix $ \mathbf{Y} $ is hermitian, has dimensions $ 3(N+1) \times 3(N+1) $, and is constructed as follows~\cite{DallAnese13,Robbins16}:
\begin{itemize}
	\item The $ |\boldsymbol{\varphi}_{lm}| \times |\boldsymbol{\varphi}_{lm}| $ off-diagonal block corresponding to the line $ (l, m) \in \mathcal{L} $ equals $ -\mathbf{Y}_{lm} \equiv - \mathbf{Z}_{lm}^{-1} $.
	\item The $ |\boldsymbol{\varphi}_{k}| \times |\boldsymbol{\varphi}_{k}| $ diagonal block corresponding to the $ k $-th bus is given by
	\begin{equation}
	[\mathbf{Y}]_{k,k} = \sum_{l \in \mathcal{N}_k} (\frac{1}{2} \mathbf{\overline{Y}}_{kl} + \mathbf{Y}_{kl})
	\end{equation}
	where $ \mathcal{N}_k := \{l : (k, l) \in \mathcal{L}\}$.
\end{itemize}

The power balance equations at node $ k \in \{\mathcal{G} \cap \mathcal{R}\} $ and phase $ \phi \in \boldsymbol{\varphi}_k $ are given by:
\begin{equation}
P_{k,\phi}^{(G)} + P_{k,\phi}^{(R)} - P_{k,\phi}^{(L)} = \real\{v_k^\phi (i_k^\phi)^\ast\}\label{pfe-1},
\end{equation}
\begin{equation}
Q_{k,\phi}^{(G)} + Q_{k,\phi}^{(R)} - Q_{k,\phi}^{(L)} = \imag\{v_k^\phi (i_k^\phi)^\ast\}\label{pfe-2}.
\end{equation}
Notice that that the balance equation for nodes without conventional generators or without RESs can be readily derived from~\eqref{pfe-1}--\eqref{pfe-2} by setting $ P_{k,\phi}^{(G)} = Q_{k,\phi}^{(G)} = 0$, or $ P_{k,\phi}^{(R)} = Q_{k,\phi}^{(R)} = 0$. Define the vectors $ \mathbf{p}_G, \mathbf{q}_G$ which collect the active and reactive powers generated by conventional generators, and let $ \mathbf{p}_R, \mathbf{q}_R$ be the vectors of  active and reactive output-powers from RESs at all nodes for  phases. Accordingly, a prototypical formulation of the AC-OPF problem for  a multi-phase power network with renewables is outlined next:
{
\begin{subequations}
	\label{OPF-F1:optim}
	\begin{align}
	& \min_{ \mathbf{v}, \mathbf{i}, \mathbf{p}_{G}, \mathbf{q}_{G}, \mathbf{p}_{R}, \mathbf{q}_{R}}
	C_g(\mathbf{p}_{G}) + C_c(\mathbf{p}_{R}) + C_i(\mathbf{q}_{R}) \label{OPF-F1:cost}\\
	&\text{subject to}\notag \\
	&\bullet\forall k \in \mathcal{N}, \phi \in \boldsymbol{\varphi}_k \notag\\
	&  \quad P_{k,\phi}^{(G)} + P_{k,\phi}^{(R)} - P_{k,\phi}^{(L)} = \real\{v_k^\phi (i_k^\phi)^\ast\} \label{OPF-F1:PBE1}\\
	&  \quad Q_{k,\phi}^{(G)} + Q_{k,\phi}^{(R)} - Q_{k,\phi}^{(L)} = \imag\{v_k^\phi (i_k^\phi)^\ast\} \label{OPF-F1:PBE2}\\
	&  \quad \underline{P}_{k,\phi}^{(G)} \leq P_{k,\phi}^{(G)} \leq \overline{P}_{k,\phi}^{(G)} \label{OPF-F1:c1}\\
	&  \quad \underline{Q}_{k,\phi}^{(G)} \leq Q_{k,\phi}^{(G)} \leq \overline{Q}_{k,\phi}^{(G)}\label{OPF-F1:c2}\\
	&  \quad |\underline{v}_{k,\phi}| \leq |v_{k,\phi}| \leq |\overline{v}_{k,\phi}| \label{OPF-F1:c3}\\
	%	&  \quad |S_{lm,\phi}| \leq |\overline{S}_{lm,\phi}| \quad \quad \forall \phi \in \boldsymbol{\varphi}_{lm},  (l,m) \in \mathcal{L} \label{OPF-F1:c4}\\
	&  \quad (P_{k,\phi}^{(R)}, Q_{k,\phi}^{(R)}) \in \mathbf{\Psi}_{k,\phi} \label{OPF-F1:c5}
	\end{align}
\end{subequations}
}
where $ \underline{P}_{k,\phi}^{(G)} $ and $ \overline{P}_{k,\phi}^{(G)} $ are the lower and upper bound on the real power generated at bus $ k $ for phase $ \phi $; $ \underline{Q}_{k,\phi}^{(G)} $ and $ \overline{Q}_{k,\phi}^{(G)} $ represents an upper and lower bounds on the reactive power injected/absorbed by a conventional generation unit at node $ k $ for phase $ \phi $; and, the constraint~\eqref{OPF-F1:c3} confine the range of the voltage magnitude of the network buses within predefined limits.
%In addition, the magnitude of apparent power flow over the line $ (l,m) $ is upper bounded by $ |\overline{S}_{lm,\phi}| $.
Notice that for buses $ \mathcal{N}\backslash \mathcal{G} $, the limits $ \underline{P}_{k,\phi}^{(G)} $, $ \overline{P}_{k,\phi}^{(G)} $, $ \underline{Q}_{k,\phi}^{(G)} $, and $ \overline{Q}_{k,\phi}^{(G)} $ are set to zero. In addition, for nodes $ \mathcal{N}\backslash \mathcal{R} $, one has $ \overline{P}_{k,\phi}^{(R)} = \overline{S}_{k,\phi}^{(R)} = 0$. The cost function~\eqref{OPF-F1:cost} is composed of three functions:
\begin{itemize}
	\item { Cost from conventional generation units:}
	\begin{equation}
	C_g(\mathbf{p}_G) = \sum_{k\in \mathcal{G}, \phi \in \boldsymbol{\varphi}_k} b_{2,k}^\phi (P_{k,\phi}^{(G)})^2 + b_{1,k}^\phi P_{k,\phi}^{(G)}
	\end{equation}
	\item { Cost of curtailment from renewables:}
	\begin{align}
	C_c(\mathbf{p}_R) & = \sum_{k\in \mathcal{R}, \phi \in \boldsymbol{\varphi}_k} c_{2,k}^\phi (\overline{P}_{k,\phi}^{(R)} - P_{k,\phi}^{(R)})^2 \nonumber \\
	& \hspace{1.8cm} + c_{1,k}^\phi (\overline{P}_{k,\phi}^{(R)} - P_{k,\phi}^{(R)})
	\end{align}
	\item {Cost of reactive support from renewables:}
	\begin{equation}
	C_i(\mathbf{q}_R) = \sum_{k\in \mathcal{R}, \phi \in \boldsymbol{\varphi}_k} d_{2,k}^\phi (Q_{k,\phi}^{(R)})^2 + d_{1,k}^\phi Q_{k,\phi}^{(R)}
	\end{equation}
\end{itemize}
Additional terms can be considered in the cost function to minimize e.g., power losses and other operational objectives.

In order to facilitate the use of the FPP-SCA algorithm, an equivalent formulation of~\eqref{OPF-F1:optim} will be introduced next. To this end, define $ \mathbf{e}_k^\phi := [\mathbf{0}_{\sum_{n=0}^{k-1} |\boldsymbol{\varphi}_{n}|}\quad \mathbf{e}_\phi\quad \mathbf{0}_{\sum_{n=k+1}^{N} |\boldsymbol{\varphi}_{n}|}] $, where $ \mathbf{e}_\phi $ is the $ \phi $-th standard canonical basis is $ \mathbb{R}^{|\boldsymbol{\varphi}_{k}|} $. Along the lines of~\cite{DallAnese13}, the following matrices are defined:
\begin{eqnarray}
\mathbf{Y}_{k,\phi} &= \frac{1}{2} (\mathbf{e}_{k,\phi} \mathbf{e}_{k,\phi}^T \mathbf{Y} + \mathbf{Y}^H \mathbf{e}_{k,\phi} \mathbf{e}_{k,\phi}^T  )\label{qmat-1},\\
\mathbf{\tilde{Y}}_{k,\phi} &= \frac{j}{2} (\mathbf{e}_{k,\phi} \mathbf{e}_{k,\phi}^T \mathbf{Y} - \mathbf{Y}^H \mathbf{e}_{k,\phi} \mathbf{e}_{k,\phi}^T)\label{qmat-2},\\
& \mathbf{M}_{k,\phi} = \mathbf{e}_{k,\phi} \mathbf{e}_{k,\phi}^T\label{qmat-3}.
\end{eqnarray}

Using~\eqref{qmat-1} and~\eqref{qmat-2}, the right hand side of the power flow equations~\eqref{pfe-1} and~\eqref{pfe-2} can be expressed as
\begin{eqnarray}
\real\{v_{k,\phi} (i_{k,\phi})^\ast\} = \mathbf{v}^H \mathbf{Y}_{k,\phi} \mathbf{v},\\
\imag\{v_{k,\phi} (i_{k,\phi})^\ast\} = \mathbf{v}^H \mathbf{\tilde{Y}}_{k,\phi} \mathbf{v},
\end{eqnarray}
%{\color{blue} The line flow matrices.}\\
while the magnitude square of the voltage phasor at bus $ k $ and phase $ \phi $ can be written in the following form:
\begin{equation}
|v_{k,\phi}|^2 = \mathbf{v}^H \mathbf{M}_{k,\phi} \mathbf{v} \, .
\end{equation}

With these definitions, the AC OPF problem~\eqref{OPF-F1:optim} can be re-written  in the following equivalent form:
{
\begin{subequations}
	\label{OPF-F2:optim}
	\begin{align}
	&	 \min_{\substack{ \mathbf{v},  \mathbf{p}_{R}, \mathbf{q}_{r}, \boldsymbol{\alpha}}}
	C_g(\boldsymbol{\alpha}) + C_c(\mathbf{p}_{R}) + C_i(\mathbf{q}_{R}) \label{OPF-F2:cost}\\
	&\text{subject to}\notag\\
	&\bullet\forall k \in \mathcal{N}, \phi \in \boldsymbol{\varphi}_k \notag\\
	&\quad \mathbf{v}^H \mathbf{Y}_{k,\phi} \mathbf{v} - P_{k,\phi}^{(R)} + P_{k,\phi}^{(L)} \leq {\alpha}_{k,\phi}  \label{OPF-F2:c1}\\
	&\quad  \underline{P}_{k,\phi}^{(G)} \leq \mathbf{v}^H \mathbf{Y}_{k,\phi} \mathbf{v} - P_{k,\phi}^{(R)} + P_{k,\phi}^{(L)}\leq \overline{P}_{k,\phi}^{(G)}   \label{OPF-F2:c2}\\
	&\quad \underline{Q}_{k,\phi}^{(G)} \leq \mathbf{v}^H \mathbf{\tilde{Y}}_{k,\phi} \mathbf{v} - Q_{k,\phi}^{(R)} + Q_{k,\phi}^{(L)}\leq \overline{Q}_{k,\phi}^{(G)} \label{OPF-F2:c3}\\
	&  \quad (|\underline{v}_{k,\phi}|)^2 \leq \mathbf{v}^H \mathbf{M}_{k,\phi} \mathbf{v} \leq (|\underline{v}_{k,\phi}|)^2 \label{OPF-F2:c4}\\
	& \quad (P_{k,\phi}^{(R)}, Q_{k,\phi}^{(R)}) \in \mathbf{\Psi}_{k,\phi}\label{OPF-F2:c5}
	%	&\bullet \forall (l,m) \in \mathcal{L}, \phi \in \boldsymbol{\varphi}_{lm} \notag\\
	%	& \quad -{t_p}_{lm,\phi} \leq \mathbf{v}^H \mathbf{Y}_{lm,\phi} \mathbf{v} \leq {t_p}_{lm,\phi}  \label{OPF-F2:c6} \\
	%	& \quad -{t_q}_{lm,\phi} \leq \mathbf{v}^H \mathbf{\tilde{Y}}_{lm,\phi} \mathbf{v} \leq {t_q}_{lm,\phi}  \label{OPF-F2:c7}\\
	%	& \quad ({t_p}_{lm,\phi})^2 + ({t_q}_{lm,\phi})^2 \leq |\overline{S}_{lm,\phi}|^2  \label{OPF-F2:c8}
	\end{align}		
\end{subequations}
}
where $ \boldsymbol{\alpha} $ is a vector that collects all $ {\alpha}_{k,\phi} $ for all $ k \in \mathcal{G} $ and $ \phi \in \boldsymbol{\varphi}_k$, {and $ {\alpha}_{k,\phi} $  represents a tight upper bound on the active power generated at node $ k $ for phase $ \phi $.}
%Also, the vectors $ \mathbf{t}_p $ and $ \mathbf{t}_q $ collect the slack variables $ t_{p_{lm,\phi}} $ and $ t_{q_{lm,\phi}} $ for $ (l, m ) \in \mathcal{L} $ and $ \phi \in \boldsymbol{\Phi}_{lm} $, respectively.
Note that, for $ k \notin \mathcal{G}$, the values of $ \underline{P}_{k,\phi}^{(G)} $ and $ \overline{P}_{k,\phi}^{(G)} $ are set to zero for all $ \phi \in \boldsymbol{\varphi}_k $. Similarly, for $ k \notin \mathcal{R}$, one has $ \overline{P}_{k,\phi}^{(R)} = \overline{S}_{k,\phi}^{(R)} = 0 $ for all $ \phi \in \boldsymbol{\varphi}_k $.
The problem~\eqref{OPF-F2:optim} is a nonconvex QCQP. Accordingly, the FPP-SCA algorithm will be utilized in the next section to to identify feasible solution of~\eqref{OPF-F2:optim} in  scenarios where existing convex relaxation-based methods may fail.

\section{Feasible Point Pursuit and Successive Convex Approximation Algorithm}
The FPP-SCA is a two-step algorithm that involves solving convex optimization problems iteratively. In the first step, we solve an inner approximation of \eqref{OPF-F2:optim} around a particular point. In order to ensure feasibility of the approximation, we add a slack variable $ s $ to the constraints and minimize $ s $ over the approximated feasible set. We then use the solution as an approximation point for the next step. If the slack variable becomes zero, we get a feasible point. In the second step, we solve a sequence of problems which are inner approximations of \eqref{OPF-F2:optim} around feasible points until convergence to a KKT point.

%\vspace{-2.5em}
\subsection{Feasible Point Pursuit}
In each iteration, the non-convex feasiblity set of \eqref{OPF-F2:optim} is replaced by a convex inner approximation. Each non-convex quadratic constraint is replaced by a convex restriction around a specific point. For instance, consider the constraint \eqref{OPF-F2:c2} which can be written as two inequalities in the following form.
\begin{subequations}
\begin{align}
\mathbf{v}^H \mathbf{Y}_{k,\phi} \mathbf{v} &\leq -P_{k,\phi}^{(L)} + P_{k,\phi}^{(R)} + \overline{P}_{k,\phi}^{(G)}\label{eqconst:C1},\\
\mathbf{v}^H (-\mathbf{Y}_{k,\phi}) \mathbf{v} &\leq P_{k,\phi}^{(L)} - P_{k,\phi}^{(R)} - \underline{P}_{k,\phi}^{(G)}\label{eqconst:C2}.
\end{align}
\end{subequations}
Both constraints are non-convex as the matrices $ \mathbf{Y}_{k}^{\phi} $ are indefinite. Consider \eqref{eqconst:C1} where the inequality can be rewritten as
\begin{equation}
 \mathbf{v}^H \mathbf{Y}_{k,\phi}^{(+)} \mathbf{v}+ \mathbf{v}^H \mathbf{Y}_{k,\phi}^{(-)} \mathbf{v} \leq -P_{k,\phi}^{(L)} + P_{k,\phi}^{(R)} + \overline{P}_{k,\phi}^{(G)}\label{quadconst}
\end{equation}
where $ \mathbf{Y}_{k,\phi}^{(+)}$ and $\mathbf{Y}_{k,\phi}^{(-)} $ are the positive semidefinite and the negative semidefinite parts of the matrix $ \mathbf{Y}_{k,\phi} $, respectively. For $ \mathbf{Y}_{k,\phi}^{(-)} $, the following inequality holds.
\begin{equation}
(\mathbf{v} - \mathbf{z})^H \mathbf{Y}_{k,\phi}^{(-)} (\mathbf{v}-\mathbf{z}) \leq 0.
\end{equation}
Then, expanding the left hand side, the following inequality can be obtained
\begin{equation}
\mathbf{v}^H \mathbf{Y}_{k,\phi}^{(-)} \mathbf{v} \leq 2 \mathbf{z}^H \mathbf{Y}_{k,\phi}^{(-)} \mathbf{v} - \mathbf{z}^H \mathbf{Y}_{k,\phi}^{(-)} \mathbf{z}.
\end{equation}
Hence, the surrogate function for the non-convex quadratic constraint \eqref{quadconst} can be defined as
\begin{equation}
 \mathbf{v}^H \mathbf{Y}_{k,\phi}^{(+)} \mathbf{v} + 2 \mathbf{z}^H \mathbf{Y}_{k,\phi}^{(-)} \mathbf{v} \leq -P_{k,\phi}^{(L)} + P_{k,\phi}^{(R)} + \overline{P}_{k,\phi}^{(G)} + \mathbf{z}^H \mathbf{Y}_{k,\phi}^{(-)} \mathbf{z} + s
\end{equation}
where the nonnegative slack variable $ s $ is added to ensure feasibility. Similarly, \eqref{eqconst:C2} is replaced by
\begin{equation}
 -\mathbf{v}^H\mathbf{Y}_{k,\phi}^{(-)} \mathbf{v} - 2 \mathbf{z}^H \mathbf{Y}_{k,\phi}^{(+)} \mathbf{v} \leq P_{k,\phi}^{(L)} - P_{k,\phi}^{(R)} - \underline{P}_{k,\phi}^{(G)} - \mathbf{z}^H \mathbf{Y}_{k,\phi}^{(+)} \mathbf{z} + s.
\end{equation}
{%\color{blue}
The problem to be solved in the $i$-th iteration can then be written as follows, where $\mathbf{z}_i$ is the optimum $\mathbf{ v }$ obtained in iteration $i-1$:
	\begin{subequations}
		\label{OPF-F3:optim}
		\begin{align}
		& \min_{\substack{ \mathbf{v}, \mathbf{p}_{R}, \mathbf{q}_{R} , s \geq 0}}
		\quad s   \label{OPF-F3:cost}\\
		& \text{subject to}\notag\\
		&\bullet{ \forall k \in \mathcal{N},\ \forall \phi \in \boldsymbol{\varphi}_k}\notag\\
		& \quad \mathbf{v}^H \mathbf{Y}_{k,\phi}^{(+)} \mathbf{v} + 2 \mathbf{z}_i^H \mathbf{Y}_{k,\phi}^{(-)} \mathbf{v} \leq\notag\\
		& \qquad -P_{k,\phi}^{(L)} + P_{k,\phi}^{(R)} + \overline{P}_{k,\phi}^{(G)} + \mathbf{z}_i^H \mathbf{Y}_{k,\phi}^{(-)} \mathbf{z}_i + s  \label{OPF-F3:c2a}\\
		& \quad \mathbf{v}^H (-\mathbf{Y}_{k,\phi}^{(-)}) \mathbf{v} - 2 \mathbf{z}_i^H \mathbf{Y}_{k,\phi}^{(+)} \mathbf{v} \leq\notag\\
		& \qquad P_{k,\phi}^{(L)} - P_{k,\phi}^{(R)} - \underline{P}_{k,\phi}^{(G)} - \mathbf{z}_i^H \mathbf{Y}_{k,\phi}^{(+)} \mathbf{z}_i + s \label{OPF-F3:c2b}\\
		& \quad \mathbf{v}^H \mathbf{\tilde{Y}}_{k,\phi}^{(+)} \mathbf{v} + 2 \mathbf{z}_i^H \mathbf{\tilde{Y}}_{k,\phi}^{(-)} \mathbf{v} \leq\notag\\
		& \qquad -Q_{k,\phi}^{(L)} + Q_{k,\phi}^{(R)} + \overline{Q}_{k,\phi}^{(G)} + \mathbf{z}_i^H \mathbf{\tilde{Y}}_{k,\phi}^{(-)} \mathbf{z}_i + s  \label{OPF-F3:c3a}\\
		& \quad \mathbf{v}^H (-\mathbf{\tilde{Y}}_{k,\phi}^{(-)}) \mathbf{v} - 2 \mathbf{z}_i^H \mathbf{\tilde{Y}}_{k,\phi}^{(+)} \mathbf{v} \leq\notag\\
		& \qquad Q_{k,\phi}^{(L)} - Q_{k,\phi}^{(R)} - \underline{Q}_{k,\phi}^{(G)} - \mathbf{z}_i^H \mathbf{\tilde{Y}}_{k,\phi}^{(+)} \mathbf{z}_i + s  \label{OPF-F3:c3b}\\
		& \quad \mathbf{v}^H \mathbf{M}_{k,\phi} \mathbf{v} \leq |\overline{v}_{k,\phi}|^2  + s\label{OPF-F3:c6a}\\	
		&  \quad 2 \mathbf{z}_i^H (-\mathbf{M}_{k,\phi}) \mathbf{v} \leq - |\underline{v}_{k,\phi}|^2 + \mathbf{z}_i^H (-\mathbf{M}_{k,\phi}) \mathbf{z}_i + s \label{OPF-F3:c6b}\\
		& \quad (P_{k,\phi}^{(R)}, Q_{k,\phi}^{(R)}) \in \mathbf{\Psi}_{k,\phi}\label{OPF-F3:c5}
%		&\bullet{ \forall (l,m) \in \mathcal{L},\ \forall \phi \in \boldsymbol{\varphi}_{lm}}\notag\\
%		& \quad \mathbf{v}^H \mathbf{Y}_{lm,\phi}^{(+)} \mathbf{v} + 2 \mathbf{z}_i^H \mathbf{Y}_{lm,\phi}^{(-)} \mathbf{v} \leq {t_p}_{lm,\phi} + \mathbf{z}_i^H \mathbf{Y}_{lm,\phi}^{(-)} \mathbf{z}_i + s  \label{OPF-F3:c7A} \\	
%		& \quad \mathbf{v}^H (-\mathbf{Y}_{lm,\phi}^{(-)}) \mathbf{v} - 2 \mathbf{z}_i^H \mathbf{Y}_{lm,\phi}^{(+)} \mathbf{v} \leq {t_p}_{lm,\phi} - \mathbf{z}_i^H \mathbf{Y}_{lm,\phi}^{(+)} \mathbf{z}_i + s \label{OPF-F3:c7B} \\
%		& \quad \mathbf{v}^H \mathbf{\tilde{Y}}_{lm,\phi}^{(+)} \mathbf{x} + 2 \mathbf{z}_i^H \mathbf{\tilde{Y}}_{lm,\phi}^{(-)} \mathbf{v} \leq {t_q}_{lm,\phi} + \mathbf{z}_i^H \mathbf{\tilde{Y}}_{lm,\phi}^{(-)} \mathbf{z}_i + s \label{OPF-F3:c8a} \\	
%		& \quad \mathbf{v}^H (-\mathbf{\tilde{Y}}_{lm,\phi}^{(-)}) \mathbf{v} - 2 \mathbf{z}_i^H \mathbf{\tilde{Y}}_{lm,\phi}^{(+)} \mathbf{v} \leq {t_q}_{lm,\phi} - \mathbf{z}^H \mathbf{\tilde{Y}}_{lm,\phi}^{(+)} \mathbf{z}_i + s \label{OPF-F3:c8b} \\
%		& \quad {t_p}_{lm,\phi}^2 + {t_q}_{lm,\phi}^2 \leq  (|\overline{S}_{lm,\phi}|)^2 \label{OPF-F3:c9}
		\end{align}	
	\end{subequations}
}	

The optimization problem \eqref{OPF-F3:optim} can be cast as SOCP which can be solved efficiently in polynomial time. Each problem instance is feasible due to the positive slack variable. This {\em feasible point pursuit} is summarized in Algorithm \ref{Algorithm: FPP-Alg}.

\SetKw{Init}{Initialization: }
\DontPrintSemicolon
\begin{algorithm}[ht]
	\Init{set $ i = 0 $}, and choose $ \mathbf{z}_0 $ to be the flat voltage profile.\;
	\Repeat{$s < \epsilon_1 $ or $ ||\mathbf{v}_i - \mathbf{v}_{i-1}|| \leq \epsilon_1 $}{
		$ \mathbf{v}_i, s  \leftarrow$ solution of \eqref{OPF-F3:optim}. \;
		$ \mathbf{z}_{i+1} \leftarrow \mathbf{v}_i $.\;
		$ i \leftarrow i+1 $.\;
	}
	\KwOut{$ \mathbf{v}_{f} \leftarrow \mathbf{v}_i $}
		\vspace{.3cm}
\caption{Feasible Point Pursuit Algorithm}
\label{Algorithm: FPP-Alg}
\end{algorithm}

%\reminder{From Nikos: but it cannot be guaranteed that $s$ will drop below $\epsilon_1$; so you need to add ``or the (relative) change in $s$ drops below $\delta_1$, or similar.}

It is clear that the value of $ s $ is nonincreasing with $ i $ as $  \mathbf{z}_i (\mathbf{v}_{i-1}) $ is always feasible while solving \eqref{OPF-F3:optim}. Despite the fact that this method in not guaranteed to find a feasible point, it always converges in the simulations to a voltage profile given by $ \mathbf{v}_f $ that is feasible. Therefore, $ \mathbf{v}_f $ is used as a starting point for the second part of our algorithm (SCA).

\subsection{Successive Convex Approximation}
Starting from a feasible point, the nonconvex feasible set is replaced at each iteration by an inner convex approximation. Similar to the FPP phase, the surrogates are formulated as convex upper bounds for the nonconvex parts of the quadratic constraints. Consequently, a monotone sequence that converges to a KKT point of the original problem \eqref{OPF-F1:optim} is generated. In each iteration, the following problem is solved
{%\color{blue}
\begin{subequations}
	\label{OPF-F4:optim}
	\begin{align}
	& \min_{\substack{ \mathbf{v}, \boldsymbol{\alpha}, \mathbf{p}_{R}, \mathbf{q}_{R}}}
	 \quad  C_g(\boldsymbol{\alpha}) + C_c(\mathbf{p}_{R}) +  C_i(\mathbf{q}_{R})  \label{OPF-F4:cost}\\
	&\text{subject to}\notag\\
	&  \quad \quad \quad \eqref{OPF-F3:c2a}-\eqref{OPF-F3:c5}\notag \quad \text{(with $s$ removed $\Leftrightarrow$ $s$ set to $0$)}\\
	&\bullet{ \forall k \in \mathcal{N},\ \forall \phi \in \boldsymbol{\varphi}_k}\notag\\
	&\mathbf{v}^H {\mathbf{Y}_{k,\phi}}^{(+)} \mathbf{v} + 2 \mathbf{z}_i^H {\mathbf{Y}_{k,\phi}}^{(-)} \mathbf{v} \leq\notag\\& \qquad \qquad -P_{k,\phi}^{(L)} + P_{k,\phi}^{(R)} + \mathbf{z}_i^H {\mathbf{Y}_{k,\phi}}^{(-)} \mathbf{z}_i + {\alpha}_{k,\phi} \label{OPF-F4:c1}
	\end{align}		
\end{subequations}
}
Note that, since the starting point is feasible, we do not add $ s $ to the surrogate constraints, or equivalently, the value of $ s $ is set to be zero. Therefore, the generated sequence is always feasible and the cost function is nonincreasing with the iterates. Algorithm \ref{Algorithm: SCA-Alg} describes the steps of the SCA phase.

\begin{algorithm}[h]
	\Init{set $ i = 0 $}, and $ \mathbf{z}_0 = \mathbf{v}_f $ .\;
	\Repeat{$\frac{\mathbf{v}_{i-1} - \mathbf{v}_{i}}{\mathbf{v}_{i-1}} < \epsilon_2 $}{
		$ \mathbf{v}_i  \leftarrow$ solution of \eqref{OPF-F4:optim}. \;
		$ \mathbf{z}_{i+1} \leftarrow \mathbf{v}_i $.\;
		$ i \leftarrow i+1 $.\;
	}
	\KwOut{$ \mathbf{v}_{opt} \leftarrow \mathbf{v}_i $}
	\vspace{.3cm}
	\caption{Successive Convex Approximation Algorithm}
		\label{Algorithm: SCA-Alg}
\end{algorithm}
%\begin{proposition}[Convergence]
%	 Let $ \{\mathbf{v}_r\} $ be a sequence generated by the FPP-SCA algorithm. Then, the whole solution sequence converges to the set $ \mathcal{K} $ that consists of all the KKT points of~\eqref{OPF-F1:optim}, i.e.,
%	 \begin{equation}
%	 \lim\limits_{r \rightarrow \inf} d^{(r)}(\mathcal{K}) \rightarrow 0
%	 \end{equation}
%	 where $ d^{(r)}(\mathcal{K}) =  \underset{\bar{\mathbf{v}}\in \mathcal{K}}{\min}\ ||\mathbf{v}_r - \bar{\mathbf{v}}||$.
%%	From ~\cite[Theorem~1]{razaviyayn-2014}, it can be shown that every limit point generated using the proposed algorithms is a KKT point. Hence, the first phase converges to a KKT point of \eqref{OPF-F3:optim}. In addition, if we start the second phase from a feasible initialization, then the sequence will converge to a KKT point of the OPF problem \eqref{OPF-F1:optim}. Note though that a KKT point of \eqref{OPF-F3:optim} is not guaranteed to be a feasible point of \eqref{OPF-F1:optim} -- in fact \cite{Mehanna-2014} contains a counter-example -- however our experience is that a feasible point is generated with high probability, if one exists (always the case in our OPF experiments).
%\end{proposition}
%\begin{proof}
%	The proof is relegated to the appendix.
%\end{proof}
\begin{claim}[Convergence] From ~\cite[Theorem~1]{razaviyayn-2014}, it can be shown that every limit point generated using the proposed algorithms is a KKT point. Hence, the first phase converges to a KKT point of \eqref{OPF-F3:optim}. In addition, if we start the second phase from a feasible initialization, then the whole sequence generated will converge to the set containing all the KKT points of the OPF problem \eqref{OPF-F1:optim}.
\end{claim}

The first part of the claim follows directly from~\cite{razaviyayn-2014}. For the second phase, if the initialization point is feasible, then the whole generated sequence will lie in the feasibility set. Because the feasible set is compact, i.e., closed and bounded, the whole converging sequence will go to the set that comprises all the KKT points of~\eqref{OPF-F1:optim}.
Note though that a KKT point of \eqref{OPF-F3:optim} is not guaranteed to be a feasible point of \eqref{OPF-F1:optim} -- in fact \cite{Mehanna-2014} contains a counter-example -- however our experience is that a feasible point is generated with high probability, if one exists (always the case in our OPF experiments).
%\reminder{From Nikos: This and the appendix have to be carefully edited, per our discussion.}

{
\subsection{Identifying Problematic Constraints}
The AC OPF problem may be infeasible under a number of operational settings, where the demand cannot be satisfied without violating voltage and/or flow constraints. When convex relaxation of the OPF problem is infeasible, it provides an infeasibility certificate for the original (nonconvex) problem. However, such relaxations typically cannot provide  informative feedback on the problematic constraints -- something valuable to the network operator to take corrective actions. Off-the-shelf solvers such as IPOPT cannot identify the problematic constraints either. }

{The FPP-SCA method~\eqref{OPF-F3:optim} seeks a feasible operating point in the first phase by minimizing the slack variable. The value of the slack variable at each iteration is in fact related to the maximum constraint violation. This method can be suitably modified to enable network operators to identify the constraints that render the overall OPF infeasible. Particularly, consider associating a slack variable with each constraint, and minimizing a cost function that is strictly increasing in the slack variables. Specifically, consider replacing problem~\eqref{OPF-F3:optim} with the following one:
	\begin{subequations}
		\label{OPF-F5:optim}
		\begin{align}
		& \min_{\substack{ \mathbf{v}, \mathbf{p}_{R}, \mathbf{q}_{R} , s \geq 0}}
		%\quad \|\boldsymbol{s}^{\overline{P}}\|_{2}^2 + \|\boldsymbol{s}^{\underline{P}}\|_{2}^2 + \|\boldsymbol{s}^{\overline{Q}}\|_{2}^2 + \|\boldsymbol{s}^{\underline{Q}}\|_{2}^2 + \|\boldsymbol{s}^{\overline{V}}\|_{2}^2 + \|\boldsymbol{s}^{\underline{V}}\|_{2}^2 \label{OPF-F5:cost}\\
		\quad \|\boldsymbol{s}\|_{2}^2 \label{OPF-F5:cost}\\
		& \text{subject to}\notag\\
		&\bullet{ \forall k \in \mathcal{N},\ \forall \phi \in \boldsymbol{\varphi}_k}\notag\\
		& \quad \mathbf{v}^H \mathbf{Y}_{k,\phi}^{(+)} \mathbf{v} + 2 \mathbf{z}_i^H \mathbf{Y}_{k,\phi}^{(-)} \mathbf{v} \leq\notag\\
		& \qquad -P_{k,\phi}^{(L)} + P_{k,\phi}^{(R)} + \overline{P}_{k,\phi}^{(G)} + \mathbf{z}_i^H \mathbf{Y}_{k,\phi}^{(-)} \mathbf{z}_i + s_{k, \phi}^{\overline{P}}  \label{OPF-F5:c2a}\\
		& \quad \mathbf{v}^H (-\mathbf{Y}_{k,\phi}^{(-)}) \mathbf{v} - 2 \mathbf{z}_i^H \mathbf{Y}_{k,\phi}^{(+)} \mathbf{v} \leq\notag\\
		& \qquad P_{k,\phi}^{(L)} - P_{k,\phi}^{(R)} - \underline{P}_{k,\phi}^{(G)} - \mathbf{z}_i^H \mathbf{Y}_{k,\phi}^{(+)} \mathbf{z}_i + s_{k, \phi}^{\underline{P}} \label{OPF-F5:c2b}\\
		& \quad \mathbf{v}^H \mathbf{\tilde{Y}}_{k,\phi}^{(+)} \mathbf{v} + 2 \mathbf{z}_i^H \mathbf{\tilde{Y}}_{k,\phi}^{(-)} \mathbf{v} \leq\notag\\
		& \qquad -Q_{k,\phi}^{(L)} + Q_{k,\phi}^{(R)} + \overline{Q}_{k,\phi}^{(G)} + \mathbf{z}_i^H \mathbf{\tilde{Y}}_{k,\phi}^{(-)} \mathbf{z}_i + s_{k, \phi}^{\overline{Q}}  \label{OPF-F5:c3a}\\
		& \quad \mathbf{v}^H (-\mathbf{\tilde{Y}}_{k,\phi}^{(-)}) \mathbf{v} - 2 \mathbf{z}_i^H \mathbf{\tilde{Y}}_{k,\phi}^{(+)} \mathbf{v} \leq\notag\\
		& \qquad Q_{k,\phi}^{(L)} - Q_{k,\phi}^{(R)} - \underline{Q}_{k,\phi}^{(G)} - \mathbf{z}_i^H \mathbf{\tilde{Y}}_{k,\phi}^{(+)} \mathbf{z}_i + s_{k, \phi}^{\underline{Q}}  \label{OPF-F5:c3b}\\
		& \quad \mathbf{v}^H \mathbf{M}_{k,\phi} \mathbf{v} \leq |\overline{v}_{k,\phi}|^2  + s_{k, \phi}^{\overline{V}}\label{OPF-F5:c6a}\\	
		&  \quad 2 \mathbf{z}_i^H (-\mathbf{M}_{k,\phi}) \mathbf{v} \leq - |\underline{v}_{k,\phi}|^2 + \mathbf{z}_i^H (-\mathbf{M}_{k,\phi}) \mathbf{z}_i + s_{k, \phi}^{\underline{V}} \label{OPF-F5:c6b}\\
		& \quad (P_{k,\phi}^{(R)}, Q_{k,\phi}^{(R)}) \in \mathbf{\Psi}_{k,\phi}\label{OPF-F5:c5}
		\end{align}	
	\end{subequations}
where $ \boldsymbol{s} $ is a vector collecting all the  slack variables. It is clear that, in this setting, the individual value of each slack relates to the violation of the respective constraint.}

{Let $ \mathbf{v}_i $, $ \mathbf{p}^{(R)}_i $, $ \mathbf{q}^{(R)}_i $ and $ \mathbf{s}_i $ denote the solution of~\eqref{OPF-F5:optim} at the $ i $-th iteration of the FPP algorithm. Then, using~\cite{razaviyayn-2014}, one can easily prove that the sequence generated by solving~\eqref{OPF-F5:optim} iteratively is convergent. When $ \mathbf{s}_i $ is all zeros at the $ i $-th iteration, the corresponding $ \boldsymbol{v}_i $ is a feasible solution for the original problem. On the other hand, if the problem is infeasible, then the slacks will converge to a non-zero vector and the positive elements of $ \mathbf{s} $ will provide a pointer to the constraints that cannot be satisfied.}

{Notice that replacing the $ 2 $-norm in the cost function~\eqref{OPF-F5:cost} by $ \|\mathbf{s}\|_\infty $ yields an  optimization problem that is equivalent to~\eqref{OPF-F3:optim}.}

\section{Test Cases and Results}
To demonstrate the efficacy of the proposed algorithm, three scenarios where {convex relaxation techniques and existing solvers for nonlinear (nonconvex) programs}  are not able to reveal feasible solutions will be considered. In the first case, we consider a single-phase equivalent model for a distribution system with high PV penetration. The ability of the proposed algorithm to minimize the curtailed power while respecting the network operational constraints will be demonstrated. Then, the three-phase model of the same distribution system will be presented. Finally, several transmission systems will be used to show the ability of the FPP-SCA algorithm to solve challenging OPF problem instances where other methods fail to find feasible voltage profiles. {In addition, the ability of the proposed algorithm to identify constraints that render  the OPF problem infeasible will be demonstrated.}

The proposed algorithm and the SDR one both employ the MATLAB-based optimization modeling package YALMIP~\cite{YALMIP} along with the interior-point solver SeDuMi~\cite{sedumi} on an Intel CPU @ $3.5$ GHz ($16$ GB RAM) computer. {For the IPOPT solver, A Julia/JuMP\footnotemark[1] Package for Power Network Optimization\footnotemark[2] was adopted to solve the single-phase OPF problems for transmission systems}. We initialize our algorithm with the flat voltage profile. In addition, we choose the values of $ \epsilon_1 $ and $ \epsilon_2 $ to be $ 10^{-11} $ and $ 10^{-5} $, respectively.

\footnotetext[1]{[Online] http://julialang.org/.}
\footnotetext[2]{A Julia/JuMP Package for Power Network Optimization. [Online] https://github.com/lanl-ansi/PowerModels.jl.}

\begin{figure}[htb]
	%\begin{center}
	\centering
	\includegraphics[scale = 0.45]{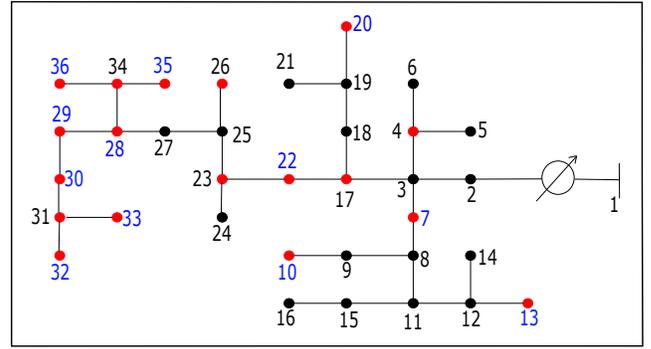}\vspace{-.5em}
	\caption{IEEE 37-node test feeder. The red nodes are the nodes with PV inverters in the single-phase model.
		 The nodes with PV units installed in the three-phase case are indexed in blue.}
	\label{fig:network1}
	%\end{center}
\end{figure}

Distribution systems with high PV penetration are likely to experience overvoltage challenges. The ability to curtail active power generated by the renewables has been shown to reliably prevent overvoltages and maintain the system operational constraints. In the first scenario, a modified version of the IEEE 37-node test feeder, shown in Fig. \ref{fig:network1}, is considered. The model is constructed by considering a single-phase equivalent feeder. Real load data measured from feeders in Anatolia, CA in August 2012~\cite{bank2013development} are used. The PV inverters are assumed to be located at the red nodes in Fig. \ref{fig:network1}, and their generation profiles are based on the real irradiance data available in~\cite{bank2013development}.

In order to show the efficacy of the proposed algorithm, two different load and irradiance profiles are considered. The first profile is taken at $ 1:00 $ PM, where the available power from PV inverters exceeds the demand. Then, the load and irradiance data at $ 7:00 $ PM in considered, where the PV inverters have very low active power. In both cases, the values of $ \overline{S}_{k,\phi}^{(R)} $ are set to be $ 2\overline{P}_{k,\phi}^{(R)} $, and the values of $ \overline{\theta}_{k,\phi} $ are set such that the minimum power factor is $ 0.7 $ for all the PV units. The limits of the voltage magnitudes $ \overline{v}_{k,\phi} $ and $ \underline{v}_{k,\phi} $ are set to be $ 1.05 $ and $ 0.95 $, respectively. Additionally, the cost function is determined by setting $ b_{2,k}^\phi = 0.1 $, $ c_{2,k}^\phi = 1 $,  $ d_{2,k}^\phi = 0.5 $, and $  b_{1,k}^\phi = c_{1,k}^\phi =  d_{1,k}^\phi = 0 $ for all $ k \in \mathcal{N} $ and $ \phi \in \boldsymbol{\varphi}_k $.

\begin{figure}[htb]
	%\begin{center}
	\centering
	\includegraphics[scale = 0.6]{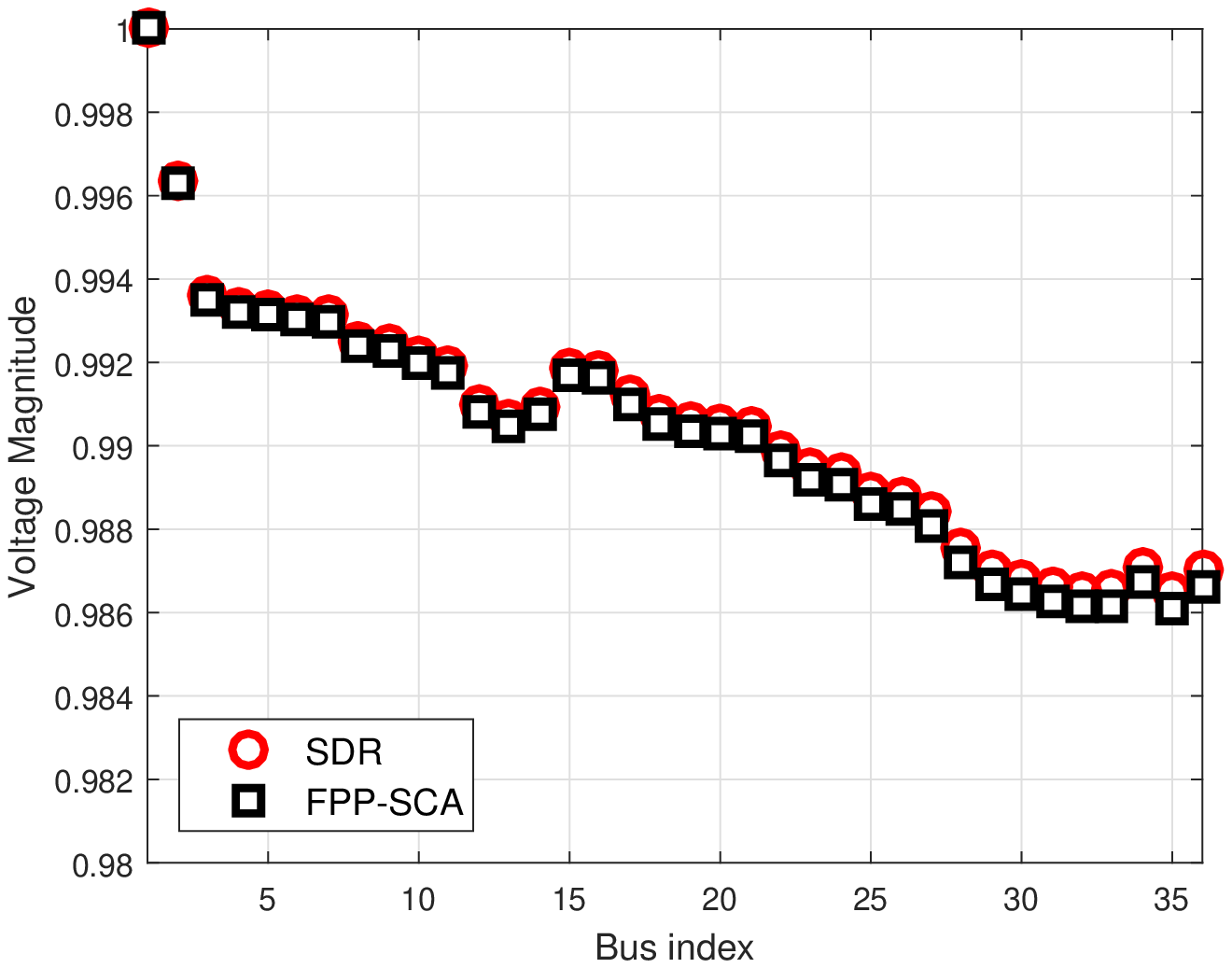}\vspace{-.5em}
	\caption{The optimal voltage profile using FPP-SCA and the SDR at $ 7:00 $ PM.	}
	\label{fig:Ph1curt0}
	%\end{center}
\end{figure}

\begin{figure}[htb]
	%\begin{center}
	\centering
	\includegraphics[scale = 0.6]{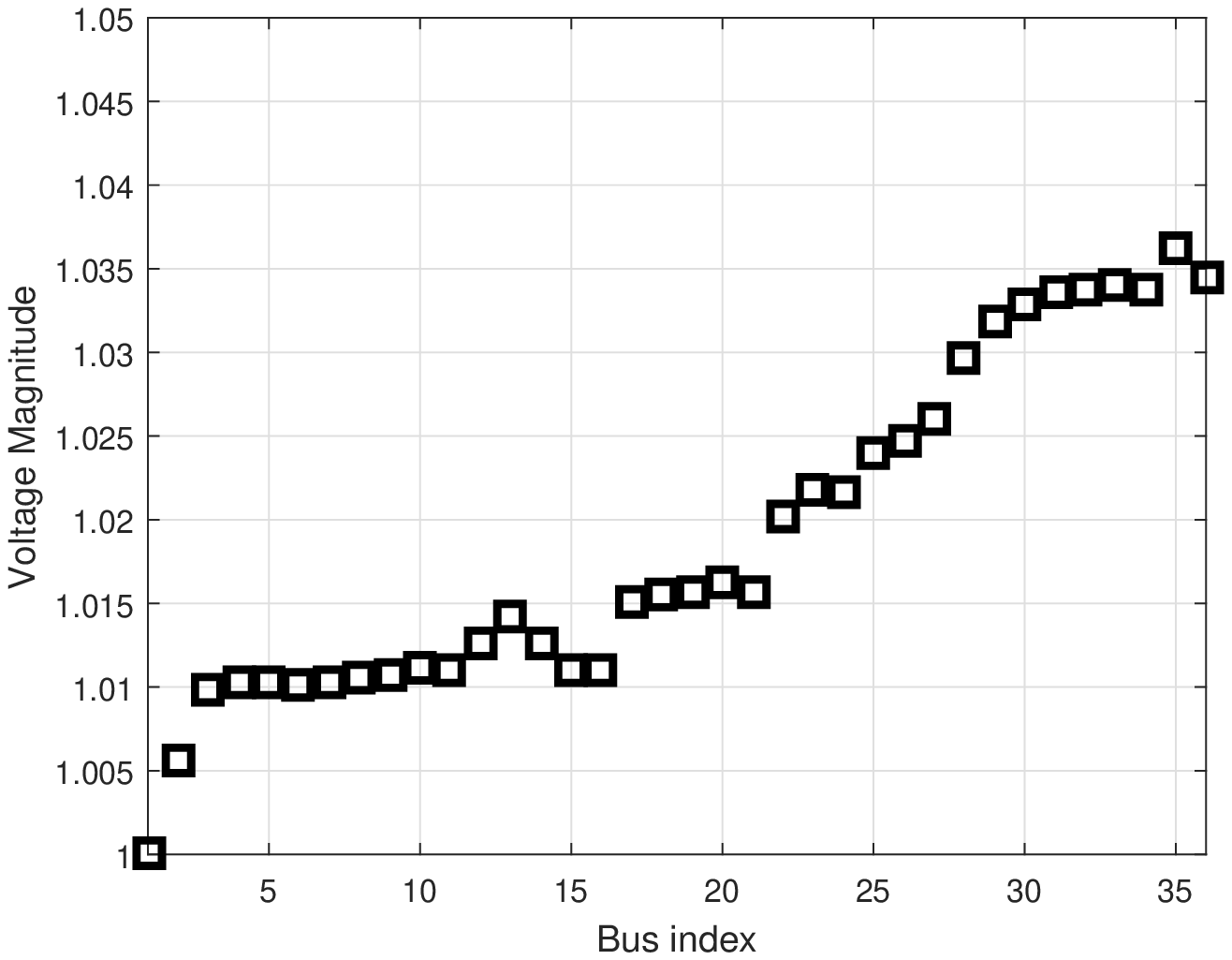}\vspace{-.5em}
	\caption{The optimal voltage profile using FPP-SCA at $ 1:00 $ PM.	}
	\label{fig:Ph1curt1}
	%\end{center}
\end{figure}

Table~\ref{tab:Results1PhCurtFeas} shows that the FPP-SCA algorithm is able to find a feasible voltage profile in both situations, while SDR is not able to find a meaningful solution when the PV penetration is high. The voltage profiles produced by the FPP-SCA and the SDR are shown in Fig. \ref{fig:Ph1curt0} in the case of low irradiance. From the depicted voltage profile, active power is drawn from node-$ 1 $ to achieve the load demand at this moment. On the other hand, the voltage profile given by the proposed algorithm at $ 1:00 $ PM is shown in Fig. \ref{fig:Ph1curt1}, where the excess of the active power generated by the PV is delivered to the transmission system connected at node-$ 1 $. Table \ref{tab:ph1curt} %\reminder{Fix table label reference}
lists the amount of the available power at each PV unit at $ 1:00 $ PM, as well as the curtailed active power resulted from the FPP-SCA solution and the injected/absorbed reactive power. Note that the power factor constraint is achieved with equality at all the PV units.

\begin{figure}
	\centering
	\begin{tabular}{|l|c|c||c|c|}
		\hline
		\multicolumn{1}{ |c|  }{\multirow{2}{*}{\textbf{Times}}}   & \multicolumn{2}{c||}{\textbf{SDR}}& \multicolumn{2}{c|}{\textbf{FPP-SCA}}\\ \cline{2-5}
		\multicolumn{1}{ |c|  }{}  & Feasibility & Cost & Feasibility & Cost\\
		\hline	
		\hline
			$ 1:00 $ PM & \texttimes     & --     & \checkmark     & $ 70486 $ \\
			$ 7:00 $ PM & \checkmark    & $ 35157 $ & \checkmark     & $ 35183 $ \\
		\hline
	\end{tabular}
	\captionof{table}{Comparison between the FPP-SCA algorithm and the SDR.}
	\label{tab:Results1PhCurtFeas}
\end{figure}

\begin{table}[htbp]
	\centering
	\begin{tabular}{cccc}
		\toprule
		%\textbf{Bus index} & \textbf{Available Power} & \textbf{Curtailed Pow}er & \textbf{Reactive Power} \\
		$ k $   & $ \overline{P}_{k,1}^{(R)} $ & $ (\overline{P}_{k,1}^{(R)} - P_{k,1}^{(R)}) $ & $ Q_{k,1}^{(R)} $\\
		\midrule
		4     & $98.60$ & $5.82$ & $ -11.46 $ \\		
		7     & $98.60$ & $5.81$ & $ -11.46 $ \\
		10    & $98.60$ & $5.78$ & $ -11.50 $ \\
		13    & $197.20$ & $ 5.71 $ &$  -11.57 $ \\
		17    & $197.20$ & $ 9.15 $ & $ -18.53  $\\
		20    & $197.20$ & $ 9.11 $ & $ -18.58 $ \\
		22    & $197.20$ & $ 12.06 $ & $ -22.92 $ \\
		23    & $197.20$ & $ 13.02 $ & $ -24.38 $ \\
		26    & $197.20$ & $ 14.53 $ & $ -26.72 $ \\
		28    & $98.60$ & $ 18.80 $ & $ -33.04  $ \\
		29    & $197.20$ & $ 18.52 $ & $-33.54 $ \\
		30    & $197.20$ & $ 18.39 $ & $ -33.81 $ \\
		31    & $ 197.20 $ & $ 18.29$ & $ -34.02 $ \\
		32    & $ 98.60 $ & $ 18.27$ & $ -34.08  $ \\
		33    & $ 197.20 $ & $ 18.26 $ & $-34.08 $ \\
		34    & $ 197.20 $ & $ 25.11 $ & $ -37.89$ \\
		35    & $ 197.20 $ & $ 41.11 $ &$  -49.40 $ \\
		36    & $ 345.10 $ & $ 25.05 $ & $ -37.93 $ \\
		\bottomrule
	\end{tabular}%
	\caption{PV inverters data for the single-phase system}
	\label{tab:ph1curt}%
\end{table}%

Next, we consider the three-phase model of the IEEE 37-node feeder. The PV units are assumed to installed at the nodes indexed in blue in Fig. \ref{fig:network1}. The PV penetration profile is adopted from the data available in~\cite{bank2013development}. An instance with high PV penetration was chosen where the SDR scheme is unable to find a feasible voltage profile. The PV penetration data is summarized in Table \ref{tab:3phcurt}, where the PV units are installed at one of the phases at selected buses. Again, we use the same constraints on $ (P_{k,\phi}^{(R)}, Q_{k,\phi}^{(R)}) $ and the cost function from the first scenario are considered.

\begin{figure}[htb]
	%\begin{center}
	\centering
	\includegraphics[scale = 0.6]{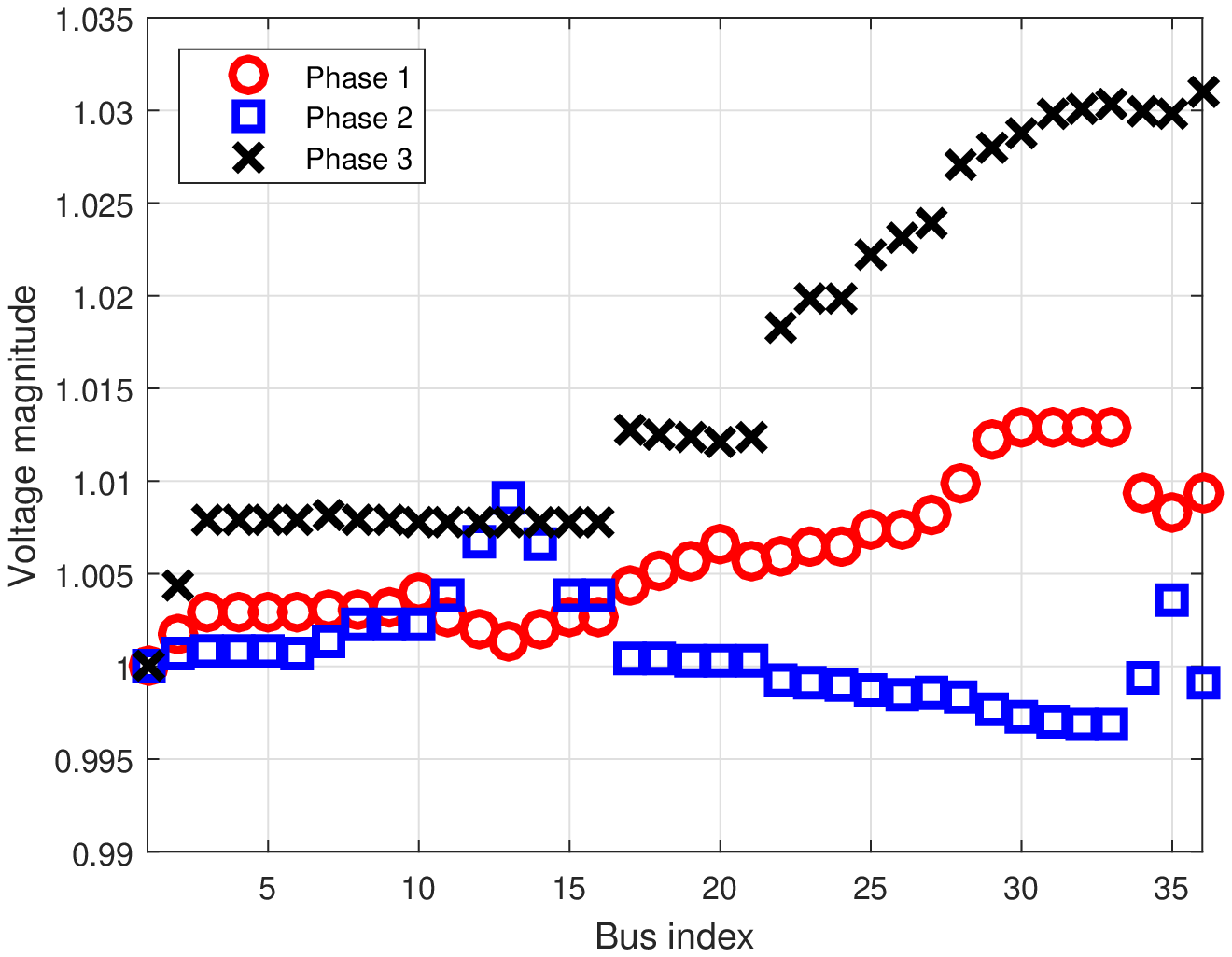}\vspace{-.5em}
	\caption{The optimal voltage profile at the three phases.	}
	\label{fig:Ph3curt1}
	%\end{center}
\end{figure}

Fig.~\ref{fig:Ph3curt1} depicts the optimal voltage profiles for the three phases across all the buses. It is clear that the voltage magnitude is high at the nodes with PV units which indicates the high power injection at these buses. Table~\ref{tab:3phcurt} lists the amount of curtailed power at the PV units, as well as the reactive power injected/absorbed by the PV inverters.

\begin{remark}
	Initializing the algorithm from the flat voltage profile in high PV penetration scenarios, the method needs about $ 1000 $ iterations in order to converge, where the subproblem can be solved in $ 5 $ and $ 2 $ seconds on average in each iteration for the single- and multi-phase systems, respectively. However, initializing the algorithm from the optimal voltage profile of close enough preceding time instance can significantly speed up the proposed algorithm. Using this strategy of warm start, the method takes only about $ 6 $ iterations (i.e., 15-30 seconds) to converge.
\end{remark}

\begin{table}[htbp]
	\centering
	\begin{tabular}{ccccc}
		\toprule
		%\textbf{Bus index}   & \textbf{Phase} & \textbf{Available Power} & \textbf{Curtailed Power} & \textbf{Reactive Power}\\
		$ k $   & $ \phi $ & $ \overline{P}_{k,\phi}^{(R)} $ & $ (\overline{P}_{k,\phi}^{(R)} - P_{k,\phi}^{(R)}) $ & $ Q_{k,\phi}^{(R)} $\\
		\midrule
		7    & $ 3 $ &  $97.86$ &  $1.23$ &  $0.32$  \\
		10    & $ 1 $ &  $97.86$ &  $0.14$ &  $-0.46$ \\
		13    & $ 2 $ &  $195.71$ &  $0$ &  $0.17$ \\
		20    & $ 1 $&  $195.71$ &  $0.23$ &  $-0.81$ \\
		22   &$ 3 $ &  $195.71$ &  $2.82$ &  $0.85$ \\
		26   & $ 3 $ &  $195.71$ &  $3.55$ &  $1.14$ \\
		28   &$ 3 $ &  $97.86$ &  $4.77$ &  $1.60$ \\
		29    & $ 1 $&  $195.71$ &  $0.02$ &  $-3.51$ \\
		30   & $ 1 $ &  $195.71$ &  $0.01$ &  $-3.87$ \\
		32    & $ 3 $&  $97.86$ &  $6.76$ &  $2.31$ \\
		33    & $ 3 $&  $195.71$ &  $6.64$ &  $2.28$\\
		35    & $ 2 $&  $195.71$ &  $0$ &  $4.82$\\
		36    & $ 3 $&  $342.5$ &  $4.78$ & $1.63$ \\
		\bottomrule
	\end{tabular}%
	\caption{PV inverters data for the three-phase system.}
	\label{tab:3phcurt}%
\end{table}%

{The ability of the proposed algorithm to solve the OPF problem instances for transmission networks is demonstrated using the test cases described in~\cite{Molzahn-2015I}. Additionally, {a modified version} of a 5-bus network presented in~\cite{bukhsh2013local} is utilized. The load and generation limits are edited to the values in Table~\ref{tab:WB5}, where the real and reactive power quantities are given in MVA and MVAr, respectively. All the other network parameters correspond to the original dataset. Table \ref{CasesTable} presents the lower bound provided by SDR, the cost of the solution produced  by FPP-SCA, and the cost obtained by IPOPT. We also compare the propose method against the moment-based relaxation~\cite{Molzahn-2015I} and the Laplacian-based approach in Table~\ref{ResultsTable}~\cite{Molzahn-2016I}.}

\begin{remark}
	For transmission networks, there are limits on the apparent power flows on the lines. Such constraints can be written as nonconvex quadratic ones after introducing slack variables. This transformation is necessary to write the OPF in QCQP form. The resulting constraints can be handled using the same way as shown before.% For more details, see~\cite{Zamzam-2016}.	
\end{remark}

\begin{table}[htbp]
	\centering
	\caption{WB5 network data.}
	\begin{tabular}{|c|cc|cccc|}
		\hline
		\multicolumn{1}{|c|}{\multirow{2}[0]{*}{Node}} & \multicolumn{2}{c|}{Load} & \multicolumn{4}{c|}{Gen. Limit} \\
		\cline{2-7}
		\multicolumn{1}{|c|}{} & $ {{P}}^{(L)} $     & $ {{Q}}^{(L)} $     & $ {\overline{P}}^{(G)} $     & $ {\underline{P}}^{(G)} $   & $ {\overline{Q}}^{(G)} $  & $ {\underline{Q}}^{(G)} $ \\
		\hline
		$ 1 $     & $ 0 $     & $ 0 $     & $ 350 $   & $ 0 $     & $ 300 $   & $ -30 $   \\
		$ 2 $     & $ 150 $   & $ 20 $    & --    & --    & --    & -- \\
		$ 3 $     & $ 150 $   & $ 20 $    & --    & --    & --    & -- \\
		$ 4 $     & $ 75 $    & $ 10 $    & --    & --    & --    & -- \\
		$ 5 $     & $ 0  $    & $ 0 $     & $ 450 $   & $ 0 $     & $ 300 $   & $ -30 $ \\
		\hline
	\end{tabular}%
	\label{tab:WB5}%
\end{table}%
\begin{figure}[htb]
	%\begin{center}
	\centering
	\includegraphics[scale = 0.8]{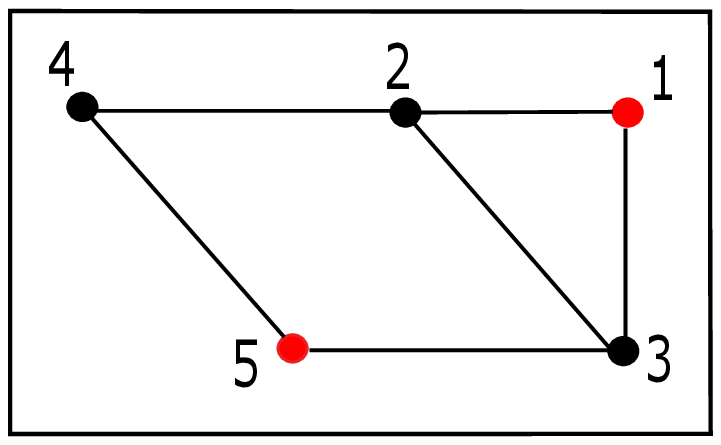}\vspace{-.5em}
	\caption{IEEE 5-node test feeder. Nodes with generators are depicted in red.}
	\label{fig:network2}
	%\end{center}
\end{figure}
%\begin{figure}
%	\centering
%	\begin{tabular}{|l|c|c|c|c|}
%		\hline
%		\textbf{Test Case}   &  \textbf{MER}& \textbf{SDR Bound} & \textbf{FPP Cost}&\textbf{IPOPT cost}\\
%		\hline		
%		\hline
%		case14Q  & $ 306.1534 $ & $3.3016\times 10^3$ & $ 3.3019\times 10^3 $& \\
%		\hline
%		case14L & $ 253.0579 $ & $9.3536\times 10^3$ &$9.3875\times 10^3$&\\
%		\hline
%		case39Q & $ 829.9366 $ & $1.0814\times 10^4$ &$1.1225\times 10^4$&\\
%		\hline
%		case39L & $ 16507.7533 $ & $4.1921\times 10^4$& $4.1974\times 10^4$&\\
%		\hline
%		case57Q & $ 423.0238 $ & $ 7.3512\times 10^3 $ &$7.3541 \times 10^3$&\\
%		\hline
%		case57L & $ 796.4583 $ & $ 4.3914\times 10^4 $ & $4.3998 \times 10^4$&\\
%		\hline
%		case118Q & $ 300.7989 $ & $ 8.1508\times 10^4 $ & $ 8.1521\times 10^4 $&\\
%		\hline
%		case118L & $ 464.9524 $ & $ 1.3391\times 10^5 $ &$ 1.3510\times 10^5 $&\\
%		\hline
%		case300 & $ 224.6997 $ & $ 7.2005\times 10^5 $ &$ 7.2016\times 10^5 $&\\
%		\hline
%	\end{tabular}
%	\captionof{table}{Test cases and results.}
%	\label{CasesTable}
%\end{figure}

\begin{figure}
	\centering
	\begin{tabular}{|l|c|c|c|}
		\hline
		\textbf{Test Case} & \textbf{SDR Bound} & \textbf{FPP Cost}&\textbf{IPOPT Cost}\\
		\hline		
		\hline
		WB5 & $ 1.1345 \times 10^3 $ & $ 1.2647 \times 10^3 $ & --\\
		\hline
		case14Q   & $3.3016\times 10^3$ & $ 3.3019\times 10^3 $& $ 3.3018\times 10^3 $\\
		\hline
		case14L  & $9.3536\times 10^3$ &$9.3875\times 10^3$&$9.3592\times 10^3$\\
		\hline
		case39Q  & $1.0814\times 10^4$ &$1.1225\times 10^4$&$1.1221\times 10^4$\\
		\hline
		case39L  & $4.1889\times 10^4$& $4.1974\times 10^4$& $4.1896\times 10^4$\\
		\hline
		case57Q  & $ 7.3472\times 10^3 $ &$7.3541 \times 10^3$&$7.3518 \times 10^3$\\
		\hline
		case57L  & $ 4.3914\times 10^4 $ & $4.3998 \times 10^4$&$4.3982 \times 10^4$\\
		\hline
		case118Q  & $ 8.1508\times 10^4 $ & $ 8.1521\times 10^4 $&$ 8.1509 \times 10^4 $\\
		\hline
		case118L  & $ 1.3391\times 10^5 $ &$ 1.3510\times 10^5 $&$ 1.3490\times 10^5 $\\
		\hline
		case300  & $ 7.1957\times 10^5 $ &$ 7.2016\times 10^5 $& $ 7.1973\times 10^5 $\\
		\hline
	\end{tabular}
	\captionof{table}{Test cases and results.}
	\label{CasesTable}
\end{figure}

{Consider the 14-,39-, 57-, 118-, and 300-bus systems (see e.g.,~\cite{Molzahn-2015I}) and {a modified version of} the 5-bus network illustrated in~\cite{bukhsh2013local}.}
These networks do not have any installed PV inverters, and hence, only traditional generation cost is considered.
{Even though the IPOPT is the most reliable software for solving the OPF problem for transmission systems, the {modified WB5 system} represents a case where IPOPT fails; on the other hand,  FPP-SCA provides a feasible (and close to optimal) solution.}
{In addition, no nonlinear solver among Trusted Region Augmented Lagrangian Multipliers (TRALM~\cite{wang2007computational}), Primal Dual Interior Point Method (PDIPM~\cite{wang2007computational}), and the Matlab Interior Point Solver (MIPS~\cite{matpower}), was able to reveal feasible solutions for all the transmission networks we tested.}
The solutions obtained using our algorithm are compared with the results of the algorithms in~\cite{Molzahn-2015I} and~\cite{Molzahn-2016I} in Table \ref{ResultsTable}. The FPP-SCA algorithm yields solutions that achieve generation costs {\em very} close to the SDR bound, in all the problem instances considered. Additionally, we compare the maximum mismatch in the nodal power injection.
We can see that the maximum mismatch in the power injection of our solution is considerably lower than the mismatch in the solutions produced by~\cite{Molzahn-2015I} and~\cite{Molzahn-2016I}. Also, whereas the solutions given by the other algorithms violate the line flow constraints by small values, the FPP-SCA algorithm is capable of finding solutions that do not violate these constraints. In these other algorithms, we may need to use higher moments to reduce the mismatch and the violation which makes the computational problem much harder. {The IPOPT solver is capable of finding solutions that are as accurate as the FPP-SCA solution; however, IPOPT may mistakenly indicate infeasibility of the OPF problem in cases where  the problem is actually feasible.}

{
\begin{figure}
	\centering
	\begin{tabular}{|l|c||c||c|}
		\hline
		\multicolumn{1}{ |c|  }{\multirow{2}{*}{\textbf{Case}}}   & \multicolumn{3}{c|}{\textbf{Maximum Injection Mismatch (MVA)}}\\ \cline{2-4}
		\multicolumn{1}{ |c|  }{}  & \textbf{MR\cite{Molzahn-2015I}} & \textbf{LA\cite{Molzahn-2016I}} & \textbf{FPP-SCA}\\
		\hline	
		\hline
		WB5 & $ 7.72 \times 10^{-9} $ & $ 3.43 $ & $ 9.07 \times 10^{-11} $\\
		\hline
		case14Q &  $ 1.08 \times 10^{-3} $ & $1.20 \times 10^{-5}$&  $ 5.15\times 10^{-8} $ \\
		\hline
		case14L & $ 5.67 \times 10^{-2} $ & $3.77 \times 10^{-5}$&  $ 2.57\times 10^{-8} $ \\
		\hline
		case39Q & $ 1.36 \times 10^{-1} $ & $ -- $& $1.26\times 10^{-4} $ \\
		\hline
		case39L & $ 4.60 \times 10^{-3} $ & $8.52 \times 10^{-3}$& $ 2.83\times 10^{-5} $ \\
		\hline
		case57Q & $ 6.49 \times 10^{-3} $ & $6.99 \times 10^{-4}$& $2.45\times 10^{-7}$ \\
		\hline
		case57L & $ 8.76 \times 10^{-4} $ & $4.42 \times 10^{-4}$& $2.37\times 10^{-6}$ \\
		\hline
		case118Q & $ 2.13 \times 10^{-1} $ & $2.98 \times 10^{-3}$& $7.52\times 10^{-6}$ \\
		\hline
		case118L & $ 4.42 \times 10^{-1} $ & $2.01 \times 10^{-3}$& $1.02\times 10^{-4}$ \\
		\hline
		case300 & $ 5.14 \times 10^{-2} $ & $7.01 \times 10^{-2}$& $7.74\times 10^{-3}$ \\
		\hline
	\end{tabular}
	\captionof{table}{Comparison between the power injection mismatch from~\cite{Molzahn-2015I} and~\cite{Molzahn-2016I} with our method.}
	\label{ResultsTable}
\end{figure}
}
%\underline{\textbf{My algorithm beats all the other algorithms.}}

%Consider the test case \textit{case39mod2} from~\cite{ArchOPF} where the Laplacian-based approach converges to an infeasible solution, i.e., the violations in the power flow limits and the nodal power injections are not acceptable. In addition, the non-convex solvers of MATPOWER package either diverge or find a local minimum solution when used to solve this case. On the other hand, our proposed algorithm converges to the best known solution for the problem~\cite{ArchOPF} which achieves a generation cost of $ 941.74\ \$/\text{hour} $. This case present a clear example of the advantage of the proposed algorithm which is able to find the solution with the lowest cost while the other methods fails to do.

\begin{figure}[htb]
	%\begin{center}
	\centering
	\includegraphics[scale = 0.6]{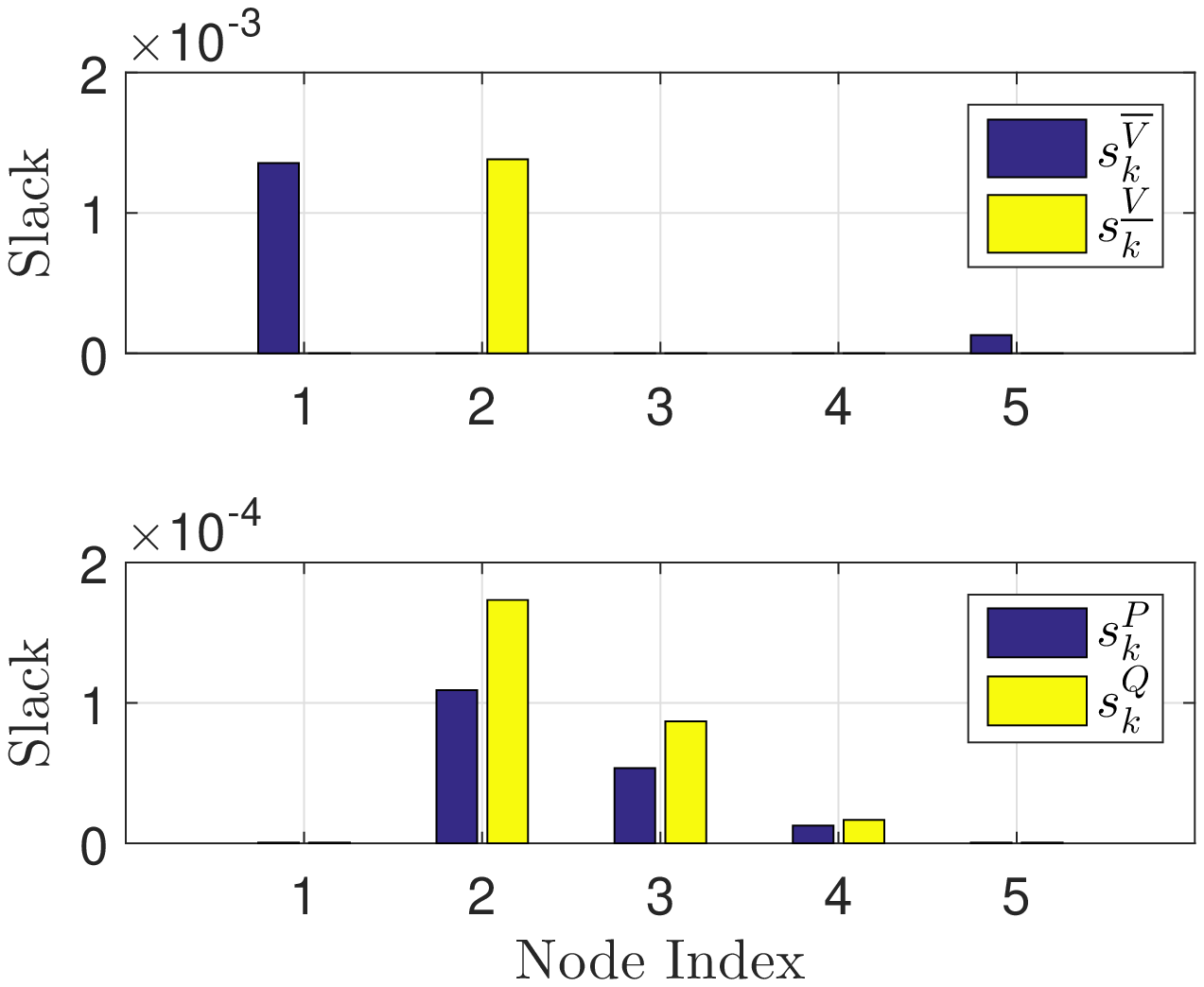}\vspace{-.5em}
	\caption{The values of slacks  of the voltage and power demand constraints for infeasible WB5.} % \textcolor{red}{Notation in the legend is different than in the rest of the paper.}}
	\label{fig:slacks}
	%\end{center}
\end{figure}

{
As an illustrative example, the networks WB5 and {\it case9mod}\footnotemark[1] are utilized next to demonstrate the ability of the proposed algorithm to identify the constraints that render the OPF infeasible. For {WB5} network, the reactive demand at node $ 2 $ is increased from  $ 20 $ MVAr to $ 70 $ MVAr. For this setting, the problem is infeasible. The value of the slack variables associated with the voltage magnitude constraints are illustrated in the upper panel of Fig.~\ref{fig:slacks}. Additionally, the slack variables associated with the loads are illustrated in the lower panel of Fig.~\ref{fig:slacks}, where $ s_k^{P}$ and $ s_k^{Q} $ are given by $\max\{s_k^{\underline{P}}, s_k^{\overline{P}}\} $ and $\max\{s_k^{\underline{Q}}, s_k^{\overline{Q}}\} $, respectively. The slack variables suggest that the upper limit of the voltage magnitude at node $ 1 $ is tight and the lower bound on the voltage magnitude of node $ 2 $ is tight. This suggests that the voltage difference between node $ 1 $ and node $ 2 $ should be larger in order to allow a higher flow of reactive power from the generator at node $ 1 $. Also, the slack variables that correspond to the load demand indicate that the demand at node $ 2 $ can not be satisfied under the existing network constraints.}

\footnotetext[1]{Available at \url{http://www.maths.ed.ac.uk/optenergy/LocalOpt/9busnetwork.html}}

\begin{figure}[htb]
	%\begin{center}
	\centering
	\includegraphics[scale = 0.6]{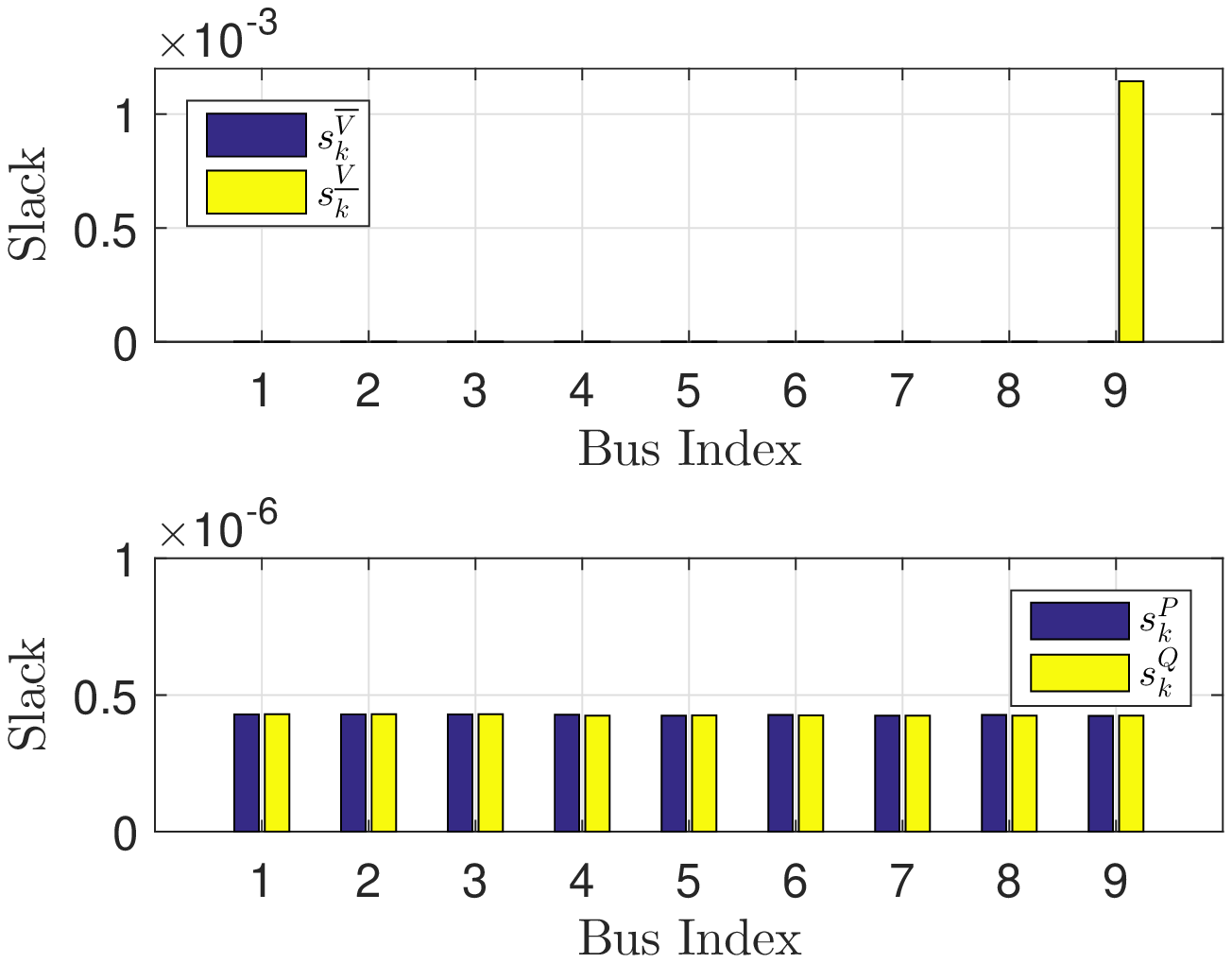}\vspace{-.5em}
	\caption{The values of the slacks of the voltage and power demand constraints for infeasible {\it case9mod}. The other bars on the top panel are not visible because their heights are $\sim 10^{-7}$.} % \textcolor{red}{Notation in the legend is different than in the rest of the paper.}}
	\label{fig:slacks2}
	%\end{center}
\end{figure}

{
A modified version of a $ 9 $-bus network~\cite{bukhsh2013local} is used next to further demonstrate the effectiveness of the FPP-SCA approach in identifying the problematic constraints. The voltage upper and lower limits were modified to be $ 1.05 $ and $ 0.95 $, respectively. In this scenario, the test case is infeasible. In Fig.~\ref{fig:slacks2}, the lower panel shows the slacks associated with the active and reactive power demand constraints. The values of these slacks are very small ($\sim 10^{-6}$), indicating that these constraints are easily satisfied. In the top panel of Fig.~\ref{fig:slacks2}, however, the slack associated with the lower limit constraint on the voltage magnitude at bus 9 is much higher, suggesting that this is the problematic constraint. Indeed, relaxing the lower limit of the voltage magnitude at bus 9 to $ 0.94 $ makes the problem feasible.
These examples represent cases where the network operator can quickly discern the source of infeasibility from the results produced by the FPP method.}

%\reminder{Add a bit more discussion / illustration for this last example? Why is it important?}

%\underline{\textbf{My algorithm is a bit slow but we want to exploit the sparsity of the quadratics.}}
%In the proposed algorithm, a convex QCQP has to be solved at each iteration which can be done optimally using several conic solvers. However, solving these problem instances is not computationally tractable for larger networks. Therefore, we can not directly apply the proposed algorithm to large test cases. Fortunately, the data matrices inherit the sparsity pattern from the underlying graph structure of the network. Thus, we can exploit this sparsity to develop an extension to the FPP-SCA algorithm that can handle these computational limitations. We are currently working on this and other speedups, and aim to deliver a scalable algorithm in the near future.

\section{Conclusions}
The AC OPF problem was considered for multi-phase networks with renewables. The problem was formulated as a nonconvex QCQP, and solved using the FPP-SCA algorithm. The proposed algorithm was shown to be effective in solving the OPF problem in many settings, including single- and three-phase models for power networks with renewables. The FPP-SCA is able to identify optimal operating points that satisfy the network constraints even under high RES penetration setups. Also, the ability of the proposed algorithm to find more accurate solutions than the moment-based relaxation and the Laplacian-based approach was demonstrated using several IEEE test cases. {Finally, the algorithm was shown to be able to identify  constraints that render the OPF problem infeasible.}

\bibliographystyle{IEEEtran}
\bibliography{IEEEabrv,OPF}

\end{document}